\documentclass[10pt,a4paper,oneside,final]{article}
\usepackage{amsmath,amsthm,amssymb,graphicx}
\newtheorem{thm}{Theorem}[section]

\newtheorem{lem}[thm]{Lemma}
\newtheorem{rem}[thm]{Remark}
\newtheorem{con}[thm]{Conjecture}
\theoremstyle{definition}
\newtheorem{defn}[thm]{Definition}
\newcommand{\Field}[1]{\mathbb{#1}}

\newcommand{\R}{\Field{R}}

\begin{document}

\title{The linking number and the writhe of uniform random walks and polygons in confined spaces}

\author{E. Panagiotou \thanks{Department of Mathematics, National Technical University of Athens, Zografou Campus, GR 15780, Athens, Greece,~\texttt{panagiotou@math.ntua.gr}},\hspace{1cm} K. C. Millett \thanks{Department of Mathematics, University of California, Santa Barbara, California 93106,~\texttt{millett@math.ucsb.edu}},\hspace{1cm} S. Lambropoulou \thanks{Department of Mathematics, National Technical University of Athens, Zografou Campus, GR 15780, Athens, Greece,~\texttt{sofia@math.ntua.gr}}
}

\maketitle
\begin{abstract}
Random walks and polygons are used to model polymers. In this paper we consider the extension of writhe, self-linking number and linking number to open chains. We then study the average writhe, self-linking and linking number of random walks and polygons over the space of configurations as a function of their length. We show that the mean squared linking number, the mean squared writhe and the mean squared self-linking number of oriented uniform random walks or polygons of length $n$, in a convex confined space, are of the form $O(n^2)$. Moreover, for a fixed simple closed curve in a convex confined space, we prove that the mean absolute value of the linking number between this curve and a uniform random walk or polygon of $n$ edges is of the form $O(\sqrt{n})$. Our numerical studies confirm those results. They also indicate that the mean absolute linking number between any two oriented uniform random walks or polygons, of $n$ edges each, is of the form $O(n)$. Equilateral random walks and polygons are used to model polymers in $\theta$-conditions. We use numerical simulations to investigate how the self-linking and linking number of equilateral random walks scale with their length.
\end{abstract}

\section{Introduction}
A polymer melt may consist of ring polymers (closed chains), linear polymers (open chains), or a mixed collection of ring and linear polymers. Polymer chains are long flexible molecules that impose spatial constraints on each other because they cannot intersect (de Gennes 1979, Rubinstein and Colby  2003). These spatial constraints, called \textit{entanglements}, affect the conformation and motion of the chains in a polymer melt and have been studied using different models of \textit{entanglement effects} in polymers (Orlandini and Whittington 2007, Tzoumanekas and Theodorou 2006).
However, a clear expression of exactly what one means by \textit{entanglement} and how one can quantify the extent of its presence remains an elusive goal.

As open polymer chains can always be moved apart, at least at long time scales, the definition of knotting or linking must refer to spatially fixed configurations, and, as a consequence, does not give rise to a topologically invariant concept. In contrast, for closed polymer chains, the concepts of knotting and linking are unchanged under continuous deformations that do not allow breakage of the chains or passage of one portion of a chain through another. To characterize the knotting of an open chain one can use the DMS method (Millett \textit{et al} 2004, Millett and Sheldon 2005) to determine the spectrum of knotting arising from the distribution of knot types created by closure of the open polygon to the ``sphere at infinity". One is then able to employ the methods of traditional knot or link theory (Kauffman 2001) to study the topology of individual constituents or of the entire collection. The knot types of individual closures or of ring polymer chains can be analyzed, for example, using the HOMFLYPT polynomial, which can distinguish the different knot and link types with high precision (Freyd \textit{et al} 1985, Przytycki and Traczyk 1987, Ewing and Millett 1997). For pairs of oriented ring polymers, for ``frozen" open polymer chains, or for a mixed frozen pair, we will see that the Gauss linking integral can be applied to calculate a linking number, which is a classical topological integer invariant in the case of pairs of ring polymers. For open or mixed pairs, the calculated quantity is a real number that is characteristic of the conformation but changes continuously under continuous deformations of the constituent chains. As polymers in solution are flexible moving molecular conformations, we will be interested in characterizations of their physical quantities at fixed moments in time as well as their time averaged properties. For the later, we will study the averages of such properties over the entire space of possible configurations.

\bigskip

By Diao \textit{et al} (1993), Diao (1995), Pippenger (1989), Sumners and Whittington (1988) we know that the probability that a polygon or open chain with $n$ edges in the cubic lattice and 3-space is unknotted goes to zero as $n$ goes to infinity. This result confirms the Frisch-Wasserman-Delbruck conjecture that long ring polymers in dilute solution will be knotted with high probability.
A stronger theorem is that the probability that any specific knot type appears as a summand of a random walk or polygon goes to 1 as $n$ goes to infinity. However the probability of forming a knot with a given type goes to zero as $n$ goes to infinity (Whittington 1992).

The knot probability of polymer molecules also depends on the extent to which the molecule is geometrically confined. This has been studied by Arsuaga \textit{et al} (2007), Tesi \textit{et al} (1994). For instance, DNA molecules confined to viral capsids have a much higher probability of being knotted. Moreover, the distribution of knot types is different from the distribution of DNA in solution (Arsuaga \textit{et al} 2002, Weber \textit{et al} 2006).

The main subject of our study concerns uniform random chains confined in a symmetric convex region in $\R^3$. Each edge of the random open chain or polygon is defined by a pair of random points in the convex region with respect to the uniform distribution. Throughout this paper we will refer to open uniform random chains and to closed uniform random chains as uniform random walks and uniform random polygons respectively. We will consider uniform random walks or polygons confined in the unit cube. The uniform random polygon (URP) model (Millett 2000) is used to  investigate the complexity of knots formed by polymer chains in confined volumes.

\bigskip

\bigskip

This paper is organized as follows. In section 2 we study the scaling of the mean squared writhe, the mean squared linking number and the mean squared self-linking number of oriented uniform random walks and polygons in confined space with respect to their length. Next, we study the scaling of the mean absolute linking number of an oriented uniform random walk or polygon and a simple closed curve, both contained in a unit cube, with respect to the length of the random walk or the polygon respectively. In section 3 we present the results of our numerical simulations which confirm the analytical results presented in section 3. Although theoretical results about the absolute linking number between two uniform random walks or polygons in confined space appear difficult to acquire, we are able to provide numerical results in section 3. Also, it is more difficult to provide analytical results for equilateral random walks or polygons than for uniform random walks or polygons. An \textit{ideal chain} is an equilateral random walk composed of freely jointed segments of equal length (equilateral) in which the individual segments have no thickness (Diao \textit{et al} 2003, Diao \textit{et al} 2005, Dobay \textit{et al} 2003) . This situation is never completely realized for real chains, but there are several types of polymeric systems with nearly ideal chains. At a special intermediate temperature, called the $\theta$-temperature, chains are in nearly ideal conformations, because the attractive and repulsive parts of monomer-monomer interactions cancel each other. In this direction, we give in section 3 numerical estimations of the scaling of the mean absolute linking number and the mean absolute self-linking number of equilateral random walks.

\section{Measures of entanglement of open curves}\label{sec}

\begin{figure}
   \begin{center}
     \includegraphics[width=0.1\textwidth]{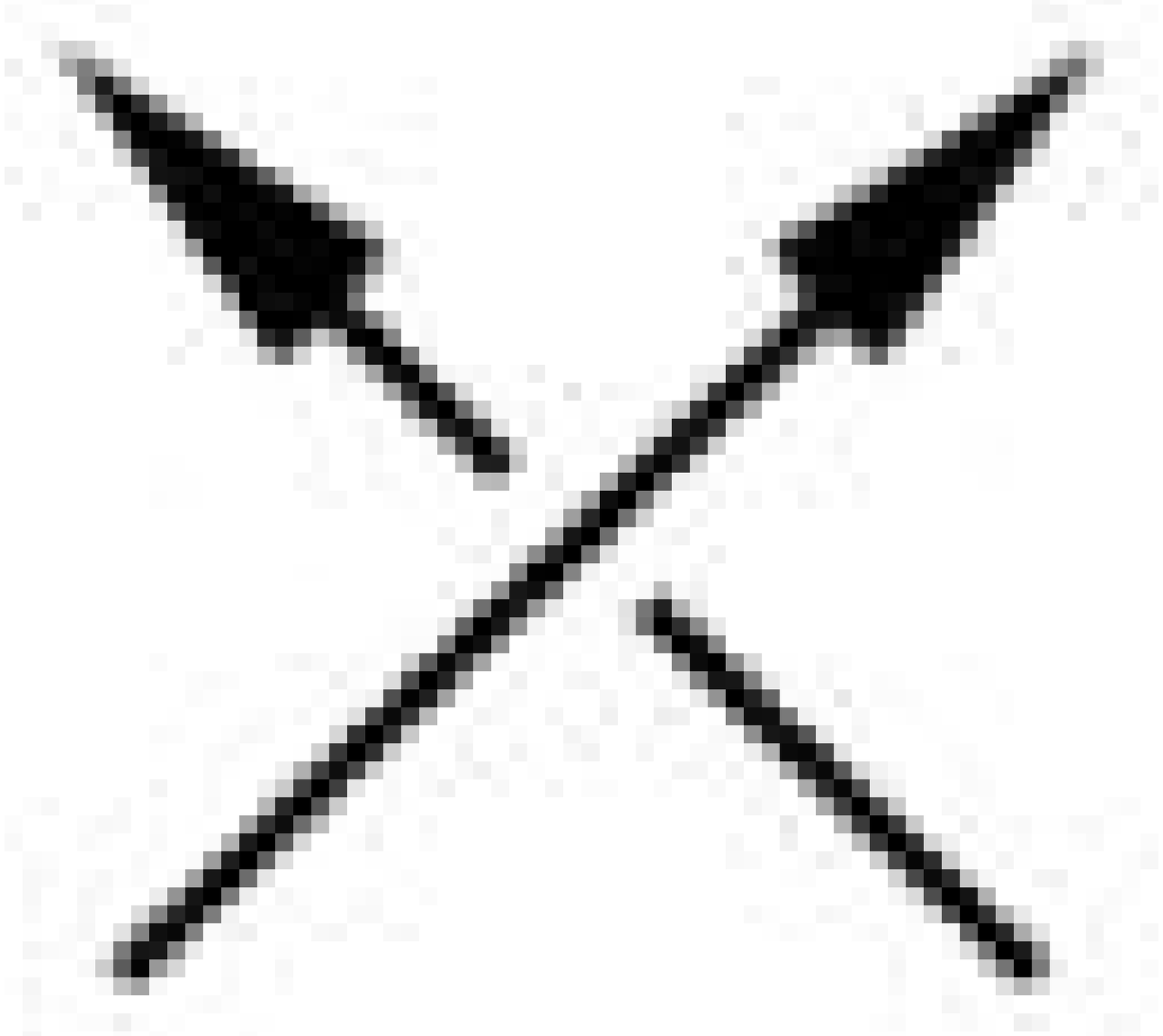}\hspace{3cm}\includegraphics[width=0.1\textwidth]{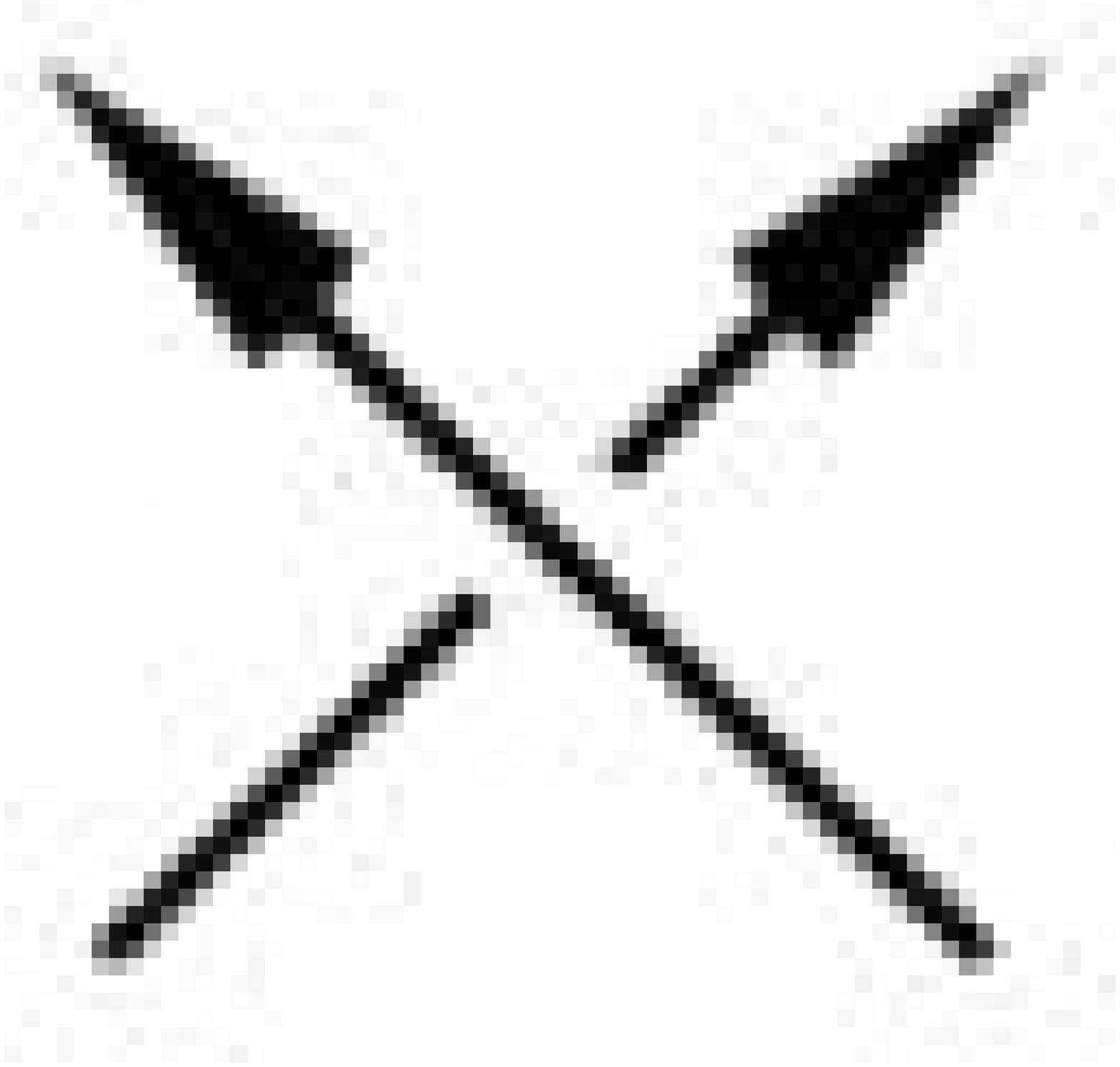}
     \caption{(a) $+1$ crossing and (b) $-1$ crossing}
     \label{sign}
   \end{center}
\end{figure}

In a generic orthogonal projection of two oriented polygonal chains, each crossing is of one of the types shown in Figure \ref{sign}. By convention, we assign $+1$ to a crossing of the first type and $-1$ to a crossing of the second type.

For a generic projection of two oriented curves $l_1,l_2$ to a plane defined by a vector $\xi\in S^2$ the \textit{linking number of a diagram}, denoted $lk_{\xi}(l_1,l_2)$, is equal to one half the algebraic sum of crossings between the projected curves.
The \textit{linking number} of two oriented curves then is equal to the average linking number of a diagram over all possible projection directions, i.e. $L(l_1,l_2)=1/4\pi(\int_{\xi\in S^2}lk_{\xi}(l_1,l_2)dS)$. This can be expressed by the Gauss linking integral for two oriented curves.

\begin{defn}\label{lk}The Gauss \textit{linking number} of two oriented curves  $l_1$ and $l_2$, whose arc-length parametrization is $\gamma_1(t),\gamma_2(s)$ respectively, is defined as a double integral over $l_1$ and $l_2$ (Gauss 1877):

\begin{equation}
L(l_1,l_2)=\frac{1}{4\pi}\int_{[0,1]}\int_{[0,1]}\frac{(\dot\gamma_1(t),\dot\gamma_2(s),\gamma_1(t)-\gamma_2(s))}{\left|\gamma_1(t)-\gamma_2(s)\right|^3}dtds
\end{equation}

where $(\dot\gamma_1(t),\dot\gamma_2(s),\gamma_1(t)-\gamma_2(s))$
is the triple product of $\dot\gamma_1(t),\dot\gamma_2(s)$ and
$\gamma_1(t)-\gamma_2(s)$.
\end{defn}

\bigskip

Similarly, for the generic orthogonal projection of one oriented curve $l$ to a plane defined by a vector $\xi\in S^2$ we define the \textit{writhe of a diagram}, denoted $Wr_{\xi}(l)$ to be equal to the algebraic sum of crossings of the projection of the curve with itself. Then the \textit{writhe} of a curve is defined as the average writhe of a diagram of the curve over all possible projections, i.e. $Wr(l)=1/4\pi(\int_{\xi\in S^2}Wr_{\xi}(l)dS)$. This can be expressed as the Gauss linking integral over one curve.

\begin{defn}\label{wr} The \textit{writhe} of an oriented curve $l$,
whose arc-length parametrization is $\gamma(t)$, is defined by the Gauss linking integral over a curve

\begin{equation}
Wr(l)=\frac{1}{2\pi}\int_{[0,1]^*}\int_{[0,1]^*}\frac{(\dot\gamma(t),\dot\gamma(s),\gamma(t)-\gamma(s))}{\left|\gamma(t)-\gamma(s)\right|^3}dtds
\end{equation}

where $[0,1]^*\times[0,1]^*=\lbrace(x,y)\in[0,1]\times[0,1]|x\neq y\rbrace$.
\end{defn}

We define the \textit{average crossing number} of a curve  $l$, whose parametrization is $\gamma(t)$, to be the average sum of crossings in a generic orthogonal projection over all possible projection directions. It is expressed by a double integral over $l$:

\begin{equation}
ACN=\frac{1}{2\pi}\int_{[0,1]}\int_{[0,1]}\frac{|(\dot\gamma(t),\dot\gamma(s),\gamma(t)-\gamma(s))|}{\left|\gamma(t)-\gamma(s)\right|^3}dtds
\end{equation}

\noindent where $(\dot\gamma(t),\dot\gamma(s),\gamma(t)-\gamma(s))$
is the triple product of $\dot\gamma(t),\dot\gamma(s)$ and
$\gamma(t)-\gamma(s)$.

\bigskip
\bigskip

We observe that the geometrical meaning of the linking number and the writhe is the same for open or closed curves. The linking number of two curves is the average over all possible projection directions of half the algebraic sum of crossings between the two components in the projection of the curves. Similarly the writhe of a curve is the average over all projection directions of the algebraic sum of crossings in the projection of the curve.

In the case of oriented closed curves the linking number is an integer topological invariant, i.e. it is invariant under isotopic moves of the curves. But in the case of oriented open curves the linking number is not topological invariant and it is not an integer. If the open curves are allowed to move continuously without intersecting each other or themselves in space, all the above measures are continuous functions in the space of configurations. Furthermore, as the endpoints of the curves move towards coincidence, the linking number, self-linking number or writhe tend to the values of those measures for the resulting closed knots or links.

\bigskip

Polymer chains are often modeled using open or closed polygonal curves. There exist several models of random walks or polygons that can be used and that are representative of the properties of different polymer melts. Since a polymer melt can take different formations in space through time, we are interested in a measure of complexity of polymer chains that will be independent of a specific configuration, and characteristic of the configuration space and how it depends on the length of the chains. In this paper, we will focus our study on uniform random walks in a confined space as this provides a simplified model for our theoretical study and will have a behaviour similar to other polymer models.

\bigskip

\section{Uniform random walks and polygons in a confined space}

The uniform random walks and polygons are modeled after the URP model, introduced by Millett (2000). In this model there are no fixed bond lengths and each coordinate of a vertex of the uniform random polygon contained in $C^3$, where $C=[0,1]$ is drawn from a uniform distribution over $[0,1]$.

The following theorem has been proved by Arsuaga \textit{et al} (2007).

\bigskip

\begin{thm}\label{uniflink}The mean squared linking number between two oriented uniform random polygons $X$ and $Y$ of $n$ edges, each contained in $C^3$, is $\frac{1}{2}n^2q$ where $q>0$. A similar result holds if $C^3$ is replaced by a symmetric convex set in $\R^3$.
\end{thm}

\bigskip

The above result is independent of the orientation of the two uniform random polygons.
Due to the weight squaring gives to larger linking number, we propose that the mean of the absolute value of the linking number between oriented uniform random walks or polygons would be a more informative measure of the expected degree of linking.

\subsection{The mean squared writhe of an oriented uniform random walk in a confined space}

In this section, we study the scaling of the writhe of an oriented uniform random walk (or polygon) contained in $C^3$.

\bigskip

We are interested in defining the average squared writhe of an $n$-step uniform random walk or polygon in confined space, where the average is taken over the entire population of open or closed uniform random walks or polygons in the confined space $C^3$. We distribute vertices according to the uniform distribution on the cube. More explicitly, the space of configurations in this case is $\Omega=[0,1]^{3(n+1)}\setminus N$, and $\Omega=[0,1]^{3n}\setminus N$ respectively, where $N$ is the set of singular configurations, i.e. when a walk or polygon intersects itself. Then $N$ is a set of measure zero (Randell 1987a, Randell 1987b).

The average writhe over the space of chains or polygons is zero as there is a sign balance occurring due to the mirror reflection involution on the space of configurations. This is why we choose to study the mean squared writhe of a uniform random walk of $n$ edges in $C^3$.

\begin{thm}\label{unifsl} The mean squared writhe of an oriented uniform random walk, or polygon of $n$ edges, each contained in $C^3$, is of the order of $O(n^2)$. Similar results hold if $C^3$ is replaced by a symmetric convex set in $\R^3$.
\end{thm}

Let us consider two (independent) oriented random edges $l_1$ and $l_2$ of an oriented uniform random polygon $P_n$ and a fixed projection plane defined by a normal vector $\xi\in S^2$.
Since the end points of the edges are independent and are uniformly distributed in $C^3$,the probability that the projections of $l_1$ and $l_2$ intersect each other is a positive number which we will call $2p$.
We define a random variable $\epsilon$ in the following way: $\epsilon=0$ if the projection of $l_1$ and $l_2$ have no intersection, $\epsilon=-1$ if the projection of $l_1$ and $l_2$ has a negative intersection and $\epsilon=1$ if the projection of $l_1$ and $l_2$ has a positive intersection. Note that, in the case the projections of $l_1$ and $l_2$ intersect, $\epsilon$ is the sign of their crossing.
Since the end points of the edges are independent and are uniformly distributed in $C^3$, we then see that $P(\epsilon=1)=P(\epsilon=-1)=p$, $E[\epsilon]=0$ and $Var(\epsilon)=E[\epsilon ^2]=2p$.

We will need the following lemma, modeled after the Lemma 1 by Arsuaga \textit{et al} (2007), concerning the case when there are four edges involved (some of them may be identical or they may have a common end point): $l_1,l_2,l_1'$ and $l_2'$. Let $\epsilon_1$ be the random number $\epsilon$ defined above between $l_1$ and $l_1'$ and let $\epsilon_2$ be the random number defined between $l_2$ and $l_2'$.

\begin{lem}\label{unilemma}
(1) If the end points of $l_1,l_2,l_1'$ and $l_2'$ are distinct,

then $E[\epsilon_1\epsilon_2]=0$(this is the case when there are eight independent random points involved).

(2) If $l_1=l_2$, and the end points of $l_1,l_1'$ and $l_2'$ are distinct(this reduces the case to where there are only three random edges with six independent points involved), then $E[\epsilon_1\epsilon_2]=0$.

(3) If $l_1$ and $l_1'$, or $l_2$ and $l_2'$ have a common end point, then $E[\epsilon_1\epsilon_2]=0$.

(4) In the case where $l_1=l_2$, the endpoints of $l_1$ and $l_1'$ and $l_1$ and $l_2'$ are distinct, and $l_1'$ and $l_2'$ share a common point (so there are only five independent random points involved in this case), let $E[\epsilon_1\epsilon_2]=u$. In the case where $l_1$ and $l_2$ share a common point,the endpoints of $l_1$ and $l_1'$ and $l_1$ and $l_2'$ are distinct, and $l_1'$ and $l_2'$ also share a common point (so there are four edges defined by six independent random points involved in this case), let $E[\epsilon_1\epsilon_2]=v$. Finally let $E[\epsilon_1\epsilon_2]=w$ in the case where $l_1$, $l_2$, $l_1'$ and $l_2'$ are consecutive (so in this case, there re four edges defined by five independent random points). Then we have $q'=3p+2(2u+v+w)>0$, where $p$ is as defined before.
\end{lem}

Note that in comparison with Lemma 1 by Arsuaga \textit{et al} (2007), we have included the case in which some of the four or three edges involved are consecutive.

\begin{proof}
(1) This is true since $\epsilon_1$ and $\epsilon_2$ are independent random variables in this case.

(2) For each configuration in which the projections of $l_1'$ and $l_2'$ both intersect the projection of $l_1$(since otherwise $\epsilon_1\epsilon_2=0$), there are eight different ways to assign the orientations to the edges. Four of them yield $\epsilon_1\epsilon_2=-1$ and four of them yield $\epsilon_1\epsilon_2=1$. Since the joint density function of the vertices involved is simply $\frac{1}{V^6}$, where $V$ is the volume of the confined space $C^3$, thus by a symmetry argument we have $E[\epsilon_1\epsilon_2]=0$.

(3) This is true since in that case $\epsilon_1=0$ or $\epsilon_2=0$.

(4) Consider the case when the polygon has six edges. Let $\epsilon_{ij}$ be the crossing sign number $\epsilon$ between the edges $l_i$ and $l_j$. Consider the variance of the summation $\sum\limits
_{1\leq i\leq 6}\sum\limits_{j>i \atop j\neq i-1,i,i+1}\epsilon_{ij}$(the summation indices are
taken modulo $6$):

\begin{equation}
\begin{split}
&V\left(\sum\limits_{1\leq i\leq6}\sum\limits_{\substack{j>i\\j\neq i-1,i,i+1}}\epsilon_{ij}\right)=E\left[\left(\sum\limits_{1\leq i\leq6}\sum\limits_{\substack{j>i\\j\neq i-1,i,i+1}}\epsilon_{ij}\right)^2\right]\\
&=\sum\limits_{1\leq i\leq6}\sum\limits_{\substack{j>i\\j\neq i-1,i,i+1}}E\left[\epsilon_{ij}^2\right]+2\sum\limits_{1\leq i\leq6}\sum\limits_{\substack{j>i\\j\neq i-2,i-1,i,i+1}}E\left[\epsilon_{ij}\epsilon_{i(j+1)}\right]\\
&+2\sum\limits_{1\leq i\leq6}\sum\limits_{\substack{j>i\\j\neq i-1,i,i+1,i+2}}E\left[\epsilon_{ij}\epsilon_{(i+1)j}\right]\\
&+2\sum\limits_{1\leq i\leq6}\sum\limits_{\substack{j>i\\j\neq i-2,i-1,i,i+1,i+2}}\left(E\left[\epsilon_{ij}\epsilon_{(i+1)(j+1)}\right]+E\left[\epsilon_{i(j+1)}\epsilon_{(i+1)j}\right]\right)\\
&+2\sum_{1\leq i\leq 6}Å\left[\epsilon_{i(i+2)}\epsilon_{(i+1)(i+3)}\right]
\end{split}
\end{equation}

Since the $\epsilon_{ij}$'s are identical random variables, i.e., they have the same distributions, each term in the first summation of the right hand side in the above yields $2p$, each term in the second summation yields $u$ (see Lemma 1), each term in the third and fourth summation yields $v$ and each term in the fifth summation yields $w$. There are $9$ terms in the first summation, there are $12$ terms in the second summation, $6$ terms in the third summation, $6$ terms in the fourth summation and $6$ terms in the fifth summation. This leads to

\begin{equation}
\begin{split}
V\left(\sum_{1\leq i\leq6}\sum_{\substack{j>i \atop j\neq i-1,i,i+1}}\epsilon_{ij}\right)&=9\cdot2p+2\left(12\cdot u+6\cdot v+6\cdot w\right)\\
&=6\left(3p+2\left(2u+v+w\right)\right).
\end{split}
\end{equation}

Since $V\left(\sum\limits_{1\leq i\leq6}\sum\limits_{\substack{j>i \atop j\neq i-1,i,i+1}}\epsilon_{ij}\right)>0$, this implies that $3p+2\left(2u+v+w\right)>0$, as claimed.
\end{proof}

\bigskip

\begin{proof}[Proof of Theorem \ref{unifsl}]

Let us consider the orthogonal projection of $P_n$ to a plane perpendicular to a random vector $\xi\in S^2$. Throughout the proof the averaging, $E$, is always over the space of configurations. The average, over the space of configurations, squared writhe of the projection of $P_n$ to that plane is given by

\begin{equation}
\begin{split}
& E\left[Wr_{\xi}^2(P_n)\right]=E\left[\left(\sum\limits_{1\leq i\leq n}\sum\limits_{\substack{j>i\\j\neq i-1,i,i+1}}\epsilon_{ij}\right)^2\right]\\
&=\sum\limits_{1\leq i\leq n}\sum\limits_{\substack{j>i\\j\neq i-1,i,i+1}}E\left[\epsilon_{ij}^2\right]+2\sum\limits_{1\leq i\leq n}\sum\limits_{\substack{j>i\\j\neq i-2,i-1,i,i+1}}E\left[\epsilon_{ij}\epsilon_{i(j+1)}\right]\\
&+2\sum\limits_{1\leq i\leq n}\sum\limits_{\substack{j>i\\j\neq i-1,i,i+1,i+2}}E\left[\epsilon_{ij}\epsilon_{(i+1)j}\right]\\
&+2\sum\limits_{1\leq i\leq n}\sum\limits_{\substack{j>i\\j\neq i-2,i-1,i,i+1,i+2}}\left(E\left[\epsilon_{ij}\epsilon_{(i+1)(j+1)}\right]+E\left[\epsilon_{i(j+1)}\epsilon_{(i+1)j}\right]\right)\\
&+2\sum_{1\leq i\leq n}Å\left[\epsilon_{i(i+2)}\epsilon_{(i+1)(i+3)}\right]\\
&=n^2(p+2(u+v))-n(3p+2(4u+5v-w))
\end{split}
\end{equation}

\noindent where $p,u,v,w$ defined as in Lemma \ref{unilemma}. By Arsuaga \textit{et al} (2007) it has been proved that $p+2(u+v)=q>0$, thus $E\left[Wr_{\xi}^2\left(P_n\right)\right]=qn^2+O\left(n\right)$. Using Lemma \ref{unilemma} we can see that $E\left[Wr_{\xi}^2\left(P_n\right)\right]$ is bounded from below by $qn^2-6qn$.

Let us now take a partition of the surface of the 2-sphere $\Delta=\lbrace I_1,I_2,\dots,I_m\rbrace$ such that the writhe of the projection of $P_n$ is constant in each $I_j, 1\leq j\leq m$.
By definition, we have that, for a sequence of partitions $\Delta_k$ such that $\mu\left(\Delta_k\right)\rightarrow0$, the mean squared writhe of $P_n$ is equal to

\begin{equation}\label{av}
\begin{split}
& E\left[Wr^2\left(P_n\right)\right]=E\left[\left(\frac{1}{4\pi}\lim_{\mu\left(\Delta_k\right)\rightarrow0}\sum_{1\leq s\leq m_k}Wr_{\xi_s}\left(P_n\right)\delta S\right)^2\right] \\
&=\frac{1}{16\pi^2}E\left[\left(\lim_{\mu\left(\Delta_k\right)\rightarrow0}\sum_{1\leq s\leq m_k} Wr_{\xi_s}\left(P_n\right)\delta S\right)^2\right] \\
&=\frac{1}{16\pi^2}E\left[\lim_{\mu\left(\Delta_k\right)\rightarrow0}\left(\sum_{1\leq s\leq m_k} Wr_{\xi_s}^2\left(P_n\right)\delta S^2+2\sum_{1\leq s,t\leq m_k}Wr_{\xi_s}\left(P_n\right)Wr_{\xi_t}\left(P_n\right)\delta S^2\right)\right] \\
&=\frac{1}{16\pi^2}\lim_{\mu\left(\Delta_k\right)\rightarrow0}E\left[\left(\sum_{1\leq s\leq m_k} Wr_{\xi_s}^2\left(P_n\right)\delta S^2+2\sum_{1\leq s,t\leq m_k}Wr_{\xi_s}\left(P_n\right)Wr_{\xi_t}\left(P_n\right)\delta S^2\right)\right] \\
& =\frac{1}{16\pi^2} \cdot \\
& \lim_{\mu\left(\Delta_k\right)\rightarrow0}\left(\sum_{1\leq s\leq m_k}E\left[Wr_{\xi_s}^2\left(P_n\right)\right]\delta S^2+2\sum_{1\leq s,t\leq m_k}E\left[Wr_{\xi_s}\left(P_n\right)\right]E\left[Wr_{\xi_t}\left(P_n\right)\right]\delta S^2\right)
\end{split}
\end{equation}

\noindent where we use the theorem of Lebesgue's dominated convergence, since the functions $S_k\left(P_n\right)=\left(\sum_{1\leq s\leq m_k}Wr_{\xi_s}\left(P_n\right)\delta S\right)^2$, $S_k:\Omega\rightarrow\R$ are measurable functions, bounded above by $\left(\sum_{1\leq s\leq m}Cr_{\xi_s}\left(P_n\right)\delta S\right)^2$, where $Cr_{\xi_s}(P_n)$ is the number of crossings of the projection of $P_n$ to the plane perpendicular to $\xi_s$ and $\left(\sum_{1\leq s\leq m_k}Cr_{\xi_s}\left(P_n\right)\delta S\right)^2\leq \left(24n^2\right)^2$.

The second term in (\ref{av}) is equal to zero, because

\begin{equation}
E\left[Wr_{\xi}\left(P_n\right)\right]=E\left[\sum_{1\leq i\leq n}\sum_{\substack{j>i \atop j\neq i-1,i,i+1}}\epsilon_{ij}\right]=\sum_{1\leq i\leq n}\sum_{\substack{j>i \atop j\neq i-1,i,i+1}}E\left[\epsilon_{ij}\right]=0.
\end{equation}

But we proved that $E\left[Wr_{\xi}^2\left(P_n\right)\right]=qn^2+O\left(n\right)$, and $E\left[Wr_{\xi}\left(P_n\right)\right]=0$,  $\forall \xi\in S^2$ so
\begin{equation}
\begin{split}
& E\left[Wr^2\left(P_n\right)\right]=\frac{1}{16\pi^2}\lim_{\mu(\Delta_k)\rightarrow0}\left(\sum_{1\leq s\leq m_k} \left(qn^2+O\left(n\right)\right)\delta S^2\right) \Rightarrow\\
&\frac{1}{16\pi^2}\left(qn^2+O\left(n\right)\right)\left(\int_{\xi\in S^2}dS\right)\leq E\left[Wr^2\left(P_n\right)\right]\leq\frac{1}{16\pi^2}\left(qn^2+O\left(n\right)\right)\left(\int_{\xi\in S^2}dS\right)^2 \\
&\Rightarrow\frac{1}{4\pi}\left(qn^2+O\left(n\right)\right)\leq E\left(Wr^2\left(P_n\right)\right)\leq qn^2+O\left(n\right)
\end{split}
\end{equation}

Note that in the case of a uniform random walk $R_n$, one has to add to the writhe the crossing between the first and last edges, otherwise the proof is the same and we obtain the same result, i.e. $E\left[Wr\left(R_n\right)^2\right]\approx qn^2+O\left(n\right)$.
\end{proof}

\subsection{The mean squared linking number of two oriented uniform random walks in confined space}

The proof of Theorem \ref{unifsl} can be easily adapted in order to provide an analysis of the rate of scaling of the mean squared linking number of open chains. Specifically, we have the following theorem, generalizing Theorem \ref{uniflink} by Arsuaga \textit{et al} (2007).

\begin{thm}\label{uniflink2}The mean squared linking number between two oriented uniform random walks $X$ and $Y$ of $n$ edges, contained in $C^3$, is is of the order of $O\left(n^2\right)$. Similar results hold if $C^3$ is replaced by a symmetric convex set in $\R^3$.
\end{thm}

\begin{proof}
For a fixed orthogonal projection of the walks to a plane perpendicular to a vector $\xi\in S^2$, adapting Theorem \ref{uniflink} of Arsuaga \textit{et al} (2007) to the case of open walks, we have $E\left[lk_{\xi}\left(X,Y\right)\right]=\frac{1}{2}n^2q+O\left(n\right)$ where $q>0$.
Then, following the proof of our Theorem \ref{unifsl}, we have that $E\left[Lk\left(X,Y\right)\right]=O\left(n^2\right)$.
\end{proof}

Note that $q=p+2\left(u+v\right)$ has the same value in all theorems. By Arsuaga \textit{et al} 2007 it has been estimated to be $q=0.0338\pm 0.024$. Our numerical results confirm this estimation.

\subsection{The mean squared self-linking number of an oriented uniform random walk or polygon}

The self-linking number was introduced to model two stranded DNA and is defined as the linking number between a curve $l$ and a translated image of that curve $l_{\epsilon}$ at a small distance $\epsilon$, i.e. $Sl\left(l\right)=L\left(l,l_{\epsilon}\right)$. This can be expressed by the Gauss integral over $\left[0,1\right]^*\times\left[0,1\right]^*=\lbrace{\left(x,y\right)\in\left[0,1\right]\times\left[0,1\right]|x\neq y\rbrace}$ by adding to it a correction term, so that it is a topological invariant of closed curves (Banchoff 1976),

\begin{eqnarray}\label{sl} SL\left(l\right)&=\frac{1}{4\pi}\int_{\left[0,1\right]^*}\int_{\left[0,1\right]^*}\frac{\left(\dot\gamma\left(t\right),\dot\gamma\left(s\right),\gamma\left(t\right)-\gamma\left(s\right)\right)}{\left|\gamma\left(t\right)-\gamma\left(s\right)\right|^3}dtds \nonumber\\
&+\frac{1}{2\pi}\int_{\left[0,1\right]}\frac{\left(\gamma'\left(t\right)\times\gamma''\left(t\right)\right)\cdot\gamma'''\left(t\right)}{|\gamma'\left(t\right)\times\gamma''\left(t\right)|^2}dt
\end{eqnarray}

The first term, in the above, is the writhe of the curve which we studied in the last section. The second term is the total torsion of the curve, $\tau\left(l\right)$ , divided by $2\pi$. This measures how much the curve deviates from being planar. The torsion of a curve can be expressed as

\begin{equation}
\tau(l)=\sum_{1\leq i\leq n}\phi_i(l)
\end{equation}

\noindent where $\phi_i(l)$ is the signed angle between the binormal vectors $B_i$ and $B_{i+1}$ defined by the edges $i-1,i,i+1$ (Banchoff 1976).

\begin{figure}
   \begin{center}
     \includegraphics[width=0.7\textwidth]{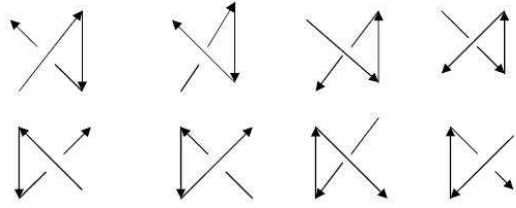}
     \caption{For all the possible configurations of the edges $i,i+1$ and $i+2$ such that $\epsilon_{i,i+2}\neq0$, we have that $\epsilon_{i,i+2}\phi_{i+1}=|\phi_{i+1}|$}
     \label{sltors}
   \end{center}
\end{figure}

The following theorem concerns the mean squared self-linking number of an oriented uniform random walk or polygon.

\begin{thm}\label{unisl} The mean squared self linking number of an oriented uniform random walk or polygon of $n$ edges, contained in $C^3$ is of the order $O(n^2)$. Similar results hold if $C^3$ is replaced by a symmetric convex set in $\R^3$.
\end{thm}

\begin{proof}
We will use the definition of the self-linking number given by (\ref{sl}), i.e. $SL(l)=Wr(l)+\frac{1}{2\pi}\sum_i\phi_i$ . The proof is based on the fact that the torsion angles $\phi_i$ $\forall i\neq1,n$ and that the products $\epsilon_{ij}\phi_i$ $\forall j\neq i+2,k\neq i+1$ are independent in the URP model.

Let $P_n$ denote a uniform random polygon in the confined space $C^3$. We project $P_n$ to a fixed plane defined by a normal vector $\xi\in S^2$.
For each pair of edges of the uniform random polygon $l_i$ and $l_j$ we define a random variable $\epsilon_{ij}$ as we did in the precious section. Then $E[\epsilon]=0$ and $E[\epsilon^2]=2p$.
For each edge we define a random variable $\phi_i$ such that $\phi_i$ is equal to the signed angle between $B_i$ and $B_{i+1}$, the normal vectors to the planes defined by the edges $i,i+1$ and $i+1,i+2$ respectively. Then $\phi_i\in[-\pi,\pi],\forall i$. Since each vertex of the uniform random polygon is chosen with respect to the uniform distribution, $\phi_i$ has equal probability of being positive or negative, thus $E(\phi_i)=0,\forall i$. Now let $E[\phi_i^2]=w$ and $E[|\phi_i|]=w'$. For each pair of edges $\phi_i,\phi_j, i\neq1,j\neq n$ we have that $E[\phi_i\phi_j]=0$, since $\phi_i,\phi_j$ are independent random variables in that case. We will now compute the mean squared self-linking number

\begin{equation}
\begin{split}
& E\left[\left(\sum_{\substack{1\leq i<j\leq n \atop j\neq i-1,i,i+1}}\epsilon_{ij}+\sum_{1\leq i\leq n}\phi_i\right)^2\right]\\
& =E\left[\left(\sum_{\substack{1\leq i<j\leq n \atop j\neq i-1,i,i+1}}\epsilon_{ij}\right)^2+\left(\sum_{1\leq i\leq n}\phi_i\right)^2+2\sum_{\substack{1\leq i<j\leq n \atop j\neq i-1,i,i+1}}\sum_{1\leq k\leq n}\epsilon_{ij}\phi_k\right]\\
& =E\left[\left(\sum_{\substack{1\leq i<j\leq n \atop j\neq i-1,i,i+1}}\epsilon_{ij}\right)^2\right]+E\left[\left(\sum_{1\leq i\leq n}\phi_i\right)^2\right]+2E\left[\sum_{\substack{1\leq i<j\leq n \atop j\neq i-1,i,i+1}}\sum_{1\leq k\leq n}\epsilon_{ij}\phi_k\right]
\end{split}
\end{equation}

\noindent In the previous section we proved that $E\left[\left(\sum\limits_{\substack{1\leq i<j\leq n \atop j\neq i-1,i,i+1}}\epsilon_{ij}\right)^2\right]=O\left(n^2\right)$.

\noindent For the second term we have that

\begin{equation}
E\left[\left(\sum_{1\leq i\leq n}\phi_i\right)^2\right]=\sum_{1\leq i\leq n}E[\phi_i^2]+2\sum_{1\leq i<j\leq n}E[\phi_i\phi_j]=w n+2E[\phi_1\phi_n]=O(n).
\end{equation}

\noindent For the third term, we proceed as follows.
If $j\neq i-2,i+2$, then $\epsilon_{ij},\phi_k$ are independent random variables for all $k$, thus $E[\epsilon_{ij}\phi_k]=0$. If $j=i+2$ then $\epsilon_{ij},\phi_k$ are independent random variables for all $k\neq i+1$ For $k=i+1$, then there are eight different cases that can occur such that $E[\epsilon_{ii+2}\phi_{i+i}]\neq0$(see Figure \ref{sltors}). All of them give $+|\phi_{i+1}|$. Since the vertices of the polygon are chosen with respect to the uniform distribution, all the cases have the same probability, thus $E[\epsilon_{ii+2}\phi_{i+1}]=E[|\phi_{i+1}|]=w'$.
So finally we have that
\begin{equation}
E[Sl^2(P_n)]=qn^2+O(n)
\end{equation}

In the case of a uniform random walk $R_n$ the self-linking number is not a topological invariant and one has to follow a similar averaging procedure as in the proof of Theorem \ref{unifsl}. Finally, the mean squared self-linking number of a uniform random walk $R_n$ is $E\left[Sl^2\left(R_n\right)\right]=O(n^2)$.
\end{proof}

\begin{rem}\rm
In the case of a uniform random polygon, Theorem \ref{unisl} can be proved using the following thinking. For a uniform random polygon $P_n$, $E\left[Sl^2\left(P_n\right)\right]=E[Lk(P_n,P_{n_{\epsilon}})^2]$, where $P_{n_{\epsilon}}$ is the translated image of $P_n$, by a distance $\epsilon$, in a direction defined by a random vector. We can then apply the same method as for the mean squared linking number of two uniform random polygons used in Theorem \ref{uniflink} by Arsuaga \textit{et al} (2007). In Lemma \ref{unilemma} by Arsuaga \textit{et al} (2007) we have to add the case of two random edges $l_i$ in $P_n$ and the translated edge $l_i'$ of $P_{n_{\epsilon}}$ such that the one is the translated image of the other. But, in that case we have $\epsilon_{ii}=0$. This does not change the rate of the scaling of the mean squared linking number and we obtain $E[Sl^2(P_n)]=O(n^2)$.

\end{rem}

\subsection{The mean absolute value of the linking number of a uniform random walk or polygon with a simple closed curve in a confined space}

In this section, following the proof of Theorem 4 by Arsuaga \textit{et al} (2007), we analyze the scaling of the absolute value of the linking number between a uniform random walk or polygon and a fixed simple closed curve in confined space.

\bigskip

\begin{thm}\label{abs2} Let $R_n$ (or $P_n$) denote an oriented uniform random walk (or polygon, respectively) of $n$ edges and $S$ a fixed simple closed curve both confined in the interior of a symmetric convex set of $\R^3$. Then the mean absolute value of the linking number between $R_n$ (or $P_n$) and $S$ has a scaling with respect to the length of the walk (or the polygon) of the form

\begin{equation}
E[\left|L(R_n,S)\right|]\approx O(\sqrt{n}).
\end{equation}

\end{thm}

\bigskip

In order to prove this theorem, we will need the following theorem from probability theory by C. Stein (Stein 1972). It is used to obtain a bound between the distribution of a sum of the terms of an $m$-dependent sequence of random variables (that is $X_1,X_2,\dots,X_s$ is independent of $X_t,X_{t+1},\dots$ provided $t-s\geq m$) and a standard normal distribution.

\begin{thm}\label{prob}
Let $x_1,x_2,\dots,x_n$ be a sequence of stationary and $m$-dependent random variables such that $E[x_i]=0,E[x_i^2]<\infty$ for each $i$ and

\begin{equation}
0< C=\lim_{n\rightarrow\infty}\frac{1}{n}E\left[\left(\sum_{1\leq i\leq n}x_i\right)^2\right]<\infty,
\end{equation}

\noindent then $\frac{1}{\sqrt{nC}}\sum_{1\leq i\leq n}x_i$ converges to the standard normal random variable. Furthermore, if we let $\Phi(a)=\frac{1}{\sqrt{2\pi}}\int_{(-\infty,a]}e^{-\frac{x^2}{2}}dx$ be the distribution function of the standard normal random variable, then we have

\begin{equation}
\left| P\left(\frac{1}{\sqrt{nC}}\sum_{1\leq i\leq n}x_i\leq a\right)-\Phi(a)\right|\leq\frac{A}{\sqrt{n}}
\end{equation}

for some constant $A>0$

\end{thm}

\begin{proof}[Proof of Theorem \ref{abs2}]
Note that the confined space can be any convex space and a simple closed curve may be of any knot type, but for simplicity we will assume that the confined space is the cube given by the set $C=\lbrace(x,y,z):-\frac{1}{2}\leq x,y,z\leq\frac{1}{2}\rbrace$ and that the simple closed curve $S$ is the circle on the $xy$-plane whose equation is $x^2+y^2=r^2$, where $r>0$ is a constant that is less than $\frac{1}{2}$.
Let $\epsilon_j$ be the sum of the $\pm1$'s assigned to the crossings between the projections of $j$th edge $l_j$ of $P_n$ and $S$, we need to take the sum since, in this case, the projections of $l_j$ may have up to two crossings with $S$. It is easy to see that $\epsilon_j=0,\pm1,\pm2$ for each $j$, the $\epsilon_j$'s have the same distributions and that, by symmetry, we have $E[\epsilon_j]=0$ for any $j$. If $|i-j|>1\bmod(n)$, then $\epsilon_i$ and $\epsilon_j$ are independent, hence we have $E[\epsilon_i\epsilon_j]=0$. Let $p'=E[\epsilon_j^2]$ and $u'=E[\epsilon_i\epsilon_{i+1}]$. Then, if $n=3$,

\begin{eqnarray}
V\left(\sum_{1\leq i\leq 3}\epsilon_{i}\right)&=E\left[\left(\sum_{1\leq i\leq 3}\epsilon_{i}\right)^2\right]=\sum_{1\leq i\leq3}E[\epsilon_i^2]+2\sum_{1\leq i,j\leq 3}E[\epsilon_i\epsilon_j]\nonumber\\
&=3p'+6u'=3(p'+2u')>0
\end{eqnarray}

Thus, we have $p'+2u'>0$, where $p'=E[\epsilon_i^2]$ and $u'=E[\epsilon_i\epsilon_j]$. It follows that
\begin{equation}\label{sum}
0<C=\frac{1}{n}E\left[\left(\sum_{1\leq j\leq n}\epsilon_j\right)^2\right]=p'+2u'
\end{equation}

for any $n$. If we ignore the last term $\epsilon_n$ in the above, then we still have

\begin{equation}
0<C=\lim_{n\rightarrow\infty}\frac{1}{n-1}E\left[\left(\sum_{1\leq j\leq n-1}\epsilon_j\right)^2\right]=p'+2u'
\end{equation}

Furthermore, the sequence $\epsilon_1,\epsilon_2,\dots,\epsilon_{n-1}$ is a stationary and 2-dependent random number sequence since the $\epsilon_j$'s have the same distributions, and what happens to $\epsilon_1,\dots,\epsilon_j$ clearly do not have any affect to what happens to $\epsilon_{j+2},\dots,\epsilon_{n-1}$ (hence they are independent).

By Theorem \ref{prob}, there exists a constant $A>0$ such that

\begin{equation}
\begin{split}
&\left|P\left(\frac{1}{\sqrt{(n-1)(p'+2u')}}\sum_{1\leq i\leq n-1}\epsilon_i\leq\alpha\right)-\Phi(\alpha)\right|\leq\frac{A}{\sqrt{n-1}}\\
&\Rightarrow\left|P\left(\sum_{1\leq i\leq n-1}\epsilon_i\leq\alpha\sqrt{(n-1)(p'+2u')}\right)-\Phi(\alpha)\right|\leq\frac{A}{\sqrt{n-1}}\\
&\Rightarrow\left|P\left(\sum_{1\leq i\leq n-1}\epsilon_i\leq w\right)-\Phi\left(\frac{w}{\sqrt{(n-1)(p'+2u')}}\right)\right|\leq\frac{A}{\sqrt{n-1}}
\end{split}
\end{equation}

where $w=\alpha\sqrt{(n-1)(p'+2u')}$.

Note that the linking number between the oriented uniform random polygon $P_n$ and $S$ is equal to the sum $\frac{1}{2}\sum_{1\leq i\leq n}\epsilon_i$.

Then as $n\rightarrow\infty$, $\frac{1}{2}\sum_{1\leq i\leq n-1}\epsilon_i\rightarrow Z$, where $Z$ is a random variable that follows the normal distribution with mean $0$ and variance $\sigma^2=\frac{1}{4}(n-1)(p'+2u')$, i.e. $N\left(0,\frac{1}{4}(n-1)(p'+2u')\right)$.
So the random variable $\left|\frac{1}{2}\sum_{1\leq i\leq n-1}\epsilon_i\right|$ follows the half normal distribution and $E\left[\left|\frac{1}{2}\sum_{1\leq i\leq n-1}\epsilon_i\right|\right]=\frac{1}{2}(2/\pi(n-1)(p'+2u'))^{1/2}=O(\sqrt{n})$.

Thus,

\begin{equation}
\left|E\left[\left|\frac{1}{2}\sum_{1\leq i\leq n-1}\epsilon_i\right|\right]-E\left[\left|\frac{1}{2}\epsilon_n\right|\right]\right|\leq E\left[\left|Lk(R_n,S)\right|\right]\leq E\left[\left|\frac{1}{2}\sum_{1\leq i\leq n-1}\epsilon_i\right|\right]+E\left[\left|\frac{1}{2}\epsilon_n\right|\right].
\end{equation}

But $E\left[\left|\frac{1}{2}\sum_{1\leq i\leq n-1}\epsilon_i\right|\right]=O(\sqrt{n})$ and $E[|\epsilon_n|]$ is a constant independent of $n$, so $E[|Lk(P_n,S)|]=O(\sqrt{n})$.

The proof carries through similarly in the case of a uniform random walk $R_n$ and a simple closed curve $S$. Note that in that case one does not have to ignore the last term in (\ref{sum}) and has to carry through an averaging procedure over all possible projections as well as in the proof of Theorem \ref{unifsl}. The result is again $E[|Lk(R_n,S)|]=O(\sqrt{n})$.
\end{proof}

\begin{rem}\rm
One can understand the above result, for the case of a uniform random polygon, and a simple closed curve using the following argument proposed by De Witt Sumners.
We know that $Lk(R_n,S)$ equals the algebraic number of times the polygon $R_n$ passes through the surface $S_1$ with $S$ as perimeter. Let $m=O(n)$ be the number of times the polygon $R_n$ passes through the surface $S_1$ with $S$ as perimeter. We can assume that $k=O(n)$ of those edges are non-consecutive (note that consecutive edges cancel each other and do not contribute to the linking number). Then we associate a variable $x_i=\pm1$ to each one of those non-consecutive edges depending upon the orientation of the edge. Observe that these are independent random variables.
Then $x_i, 1\leq i\leq k$ is a one-dimensional random walk, thus the distance of the starting point of the random walk of $k=O(n)$ steps is of the order $O(\sqrt{n})$, thus $|Lk(R_n,S)|=O(\sqrt{n})$.
\end{rem}

\bigskip

\begin{rem}\rm
We stress that it is the mean absolute value of the linking number of two uniform random walks or polygons in confined space and the mean absolute writhe and self-linking number of a uniform random walk or polygon in confined space that are of greatest interest for us. Although their analysis is much more challenging, we propose that these provide a clearer picture of the scaling of the quantities associated with the geometry and topology of the chains or polygons. We study these numerically in the next section. It would be interesting for future work to prove the numerical scaling analytically.
\end{rem}

\bigskip

\section{Numerical Results}

\bigskip

In this section we will describe results obtained by simulations of uniform random walks and polygons in confined space and of equilateral random walks. First we consider the scaling of the mean squared writhe, $E[Wr^2]$, and the mean absolute value of the writhe, $E\left[\left|Wr\right|\right]$ of an oriented  uniform random walk and polygon of $n$ edges in confined space. Then we study the scaling of the mean absolute value of the linking number $E[|Lk|]$ between an oriented uniform random walk or polygon of $n$ edges and a fixed oriented simple closed curve; and the mean absolute value of the linking number between two oriented uniform random walks or polygons of $n$ edges each. Finally we study the scaling of the mean absolute value of the linking number between two oriented equilateral random walks of $n$ edges whose starting points coincide, $\langle ALN\rangle$, and the scaling of the mean absolute value of the self-linking number of an oriented equilateral random walk of $n$ edges, with respect to the number of edges, $\langle ASL\rangle$.

\subsection{Generation of data}

To generate uniform random walks and polygons confined in $C^3$, each coordinate of a vertex of the uniform random walk was drawn from a uniform distribution on $[0,1]$, and to generate equilateral random walks, each edge vector was drawn from a uniform distribution on $S^2$.

For the computation of the linking number or the writhe of uniform or equilateral random walks or polygons, we used the algorithm by Klenin and Langowski (2000), which is based on the Gauss integral. For each pair of edges $e_1,e_2$, their linking number
is computed as the signed area of two antipodal quadrangles defined by the two edges over the area of the $2$-sphere.

We estimated the linking, the writhe and the self-linking numbers between oriented uniform random walks and polygons by analyzing pairs of 10 subcollections of 500 oriented uniform random walks or polygons ranging from 10 edges to 100 edges by a step size of 10 edges,
for which we calculated the mean and then computed the mean of the 10 means for our estimate. We did the same for the study of the linking number and the self-linking number of equilateral random walks.
For the computation of the scaling of the linking number between an oriented uniform random walk or polygon and an oriented simple closed planar curve, first we considered a fixed square and we analyzed 10 subcollections of 500 oriented uniform random walks or polygons ranging from 10 edges to 100 edges by a step size of 10 edges. In order to illustrate that the result holds for any fixed knot, we also considered a fixed trefoil and we did the same analysis.

\subsection{Analysis of data}

In this section first we analyze our data on the mean squared writhe of an oriented uniform random walk and for an oriented uniform random polygon in confined space. We fit our data to a function of the form predicted in Theorem \ref{unifsl} using the method of least-squares. Then we check the variance and the coefficient of determination. The coefficient of determination takes values between 0 and 1 and is an indicator of how well the curve fits the data. We continue further our study of the writhe of uniform random walks and polygons by analyzing our data on the mean absolute value of the writhe of a uniform random walk or polygon in confined space. We do not have an analytical result for this scaling, but the data are very well fitted to a linear function. Thus, we observe that the distribution of the writhe of a uniform random walk or polygon has the property that the scaling of the mean squared writhe is equal to the scaling of the squared mean absolute writhe. We know that if a random variable $X$ follows the normal distribution, $X\sim N(0,\sigma^2)$, then $|X|$ follows the half-normal distribution and $E[|X|]=\sigma\sqrt{2/\pi}$. Then $|X|/\sigma$ follows the $\chi$-distribution and $X^2/\sigma^2$ follows the $\chi^2$-distribution with one degree of freedom and mean $1$. Thus, $E[X^2/\sigma^2]=1$, $E[X^2]=\sigma^2$, hence we have that $E[X^2]\approx E[|X|]^2$. This indicates that the writhe of a uniform random walk or polygon in confined space follows the normal distribution.

By Arsuaga \textit{et al} (2007) it has been proved and confirmed numerically that the mean squared linking number of two oriented uniform random polygons in confined space has a scaling of the form $O(n^2)$. We stress that in the proof of this theorem we have used the fact that the linking number of two closed oriented curves is independent of the projection. For two oriented uniform random walks (open curves) in confined space we have proved that the mean squared linking number has a scaling of the form $O(n^2)$ (Theorem \ref{uniflink}). This is confirmed by our data. Since we are interested in the mean absolute value of the linking between two uniform random walks or polygons in confined space but do not have an analytical result study is challenging, we study this numerically. Our data are well fitted to a linear function. This matches the observation of Arsuaga \textit{et al} (2007) that the linking number follows the normal distribution.

In the same subsection, we study the scaling of the mean absolute linking number of a uniform random walk or polygon and a fixed simple closed curve in confined space. In Theorem \ref{abs2} we proved that it has a scaling of the form $O(\sqrt{n})$ and this is confirmed by our data. We note that this is different from the mean absolute linking number of two uniform random walks or polygons in confined space for which our numerical results indicate a scaling of the form $O(n)$.

Next we analyze numerical data for the mean squared self-linking number of uniform random walks and polygons. The data strongly support a scaling of the form $O(n^2)$ and, thus, confirm Theorem \ref{unisl}. We analyze our data for the mean absolute value of the self-linking number of a uniform random walk or polygon in confined space. Again, we observe a scaling of the form $O(n)$ which suggests that the self-linking number of a uniform random walk or polygon in confined space follows the normal distribution.

We conclude that all the above measures follow the same type of distribution. Furthermore this distribution has the property that the square root of the mean squared random variable is equal to the mean absolute random variable. This strengthens our intuition that the linking number follows a normal distribution.

\subsubsection{Mean squared writhe and mean absolute writhe of an oriented uniform random walk or polygon in a confined space}

\bigskip

Our first numerical study concerns the writhe of a uniform random walk or polygon. By Theorem \ref{unifsl} the mean squared writhe of an oriented uniform random polygon grows at a rate $E[Wr^2]\approx qn^2+O(n)$. For comparison with this analytical result we calculated the mean squared writhe of an oriented uniform random walk of varying length and of an oriented uniform random polygon of varying length.

\begin{figure}
   \begin{center}
     \includegraphics{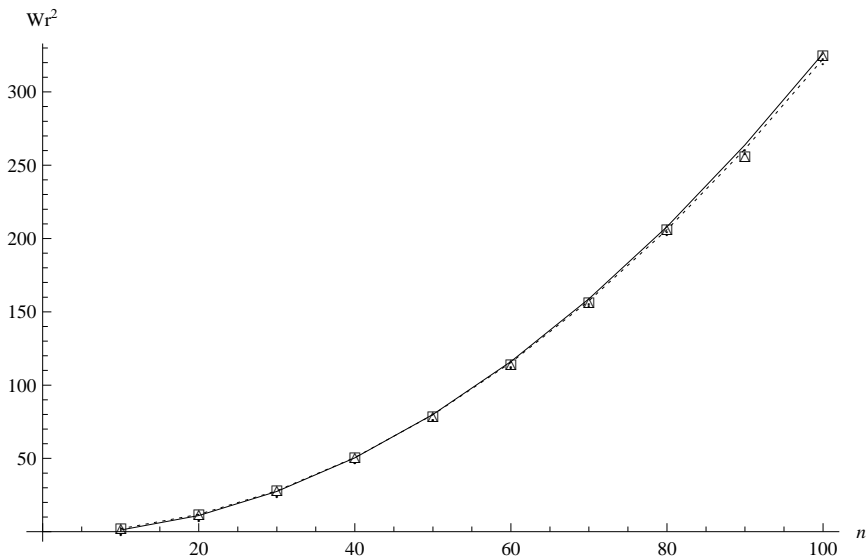}
     \caption{The mean squared writhe of uniform random walks and polygons. Values obtained by computer simulations are shown by triangles and squares respectively. The black curve is the graph of the mean squared writhe of a uniform random polygon and the dashed curve is the graph of the mean squared writhe of a uniform random walk (open chain).}
     \label{meansqwr}
   \end{center}
\end{figure}

Results are shown in Figure \ref{meansqwr}. The black curve in the figure illustrates the scaling of the mean squared writhe of a uniform random polygon with respect to its number of edges and is fitted to a function of the form $qn^2+a$ where $q$ is estimated to be $0.0329\pm0.0002$ and $a$ is estimated to be $-2.1293\pm0.9466$, with a coefficient of determination $R^2=0.9998$.
Thus the estimate given in the theorem is strongly supported by the data.
The dashed curve in the figure illustrates the scaling of the mean squared writhe of a uniform random walk with respect to its number of edges and is fitted to a function of the form $qn^2+a$ where $q$ is estimated to be $0.0324\pm0.0002$ and $a$ is estimated to be $-1.1499\pm0.9948$, with a coefficient of determination $R^2=0.9997$. Hence, the estimate given in the theorem is strongly supported by the data.

Note that $q$ was estimated to be equal to $0.0329\pm0.0002$ and $0.0324\pm0.0002$ respectively, which coincides with the predicted value of $q$ given by Arsuaga \textit{et al} (2007), i.e. $q=0.0338\pm0.024$.

\bigskip

\begin{figure}
   \begin{center}
     \includegraphics{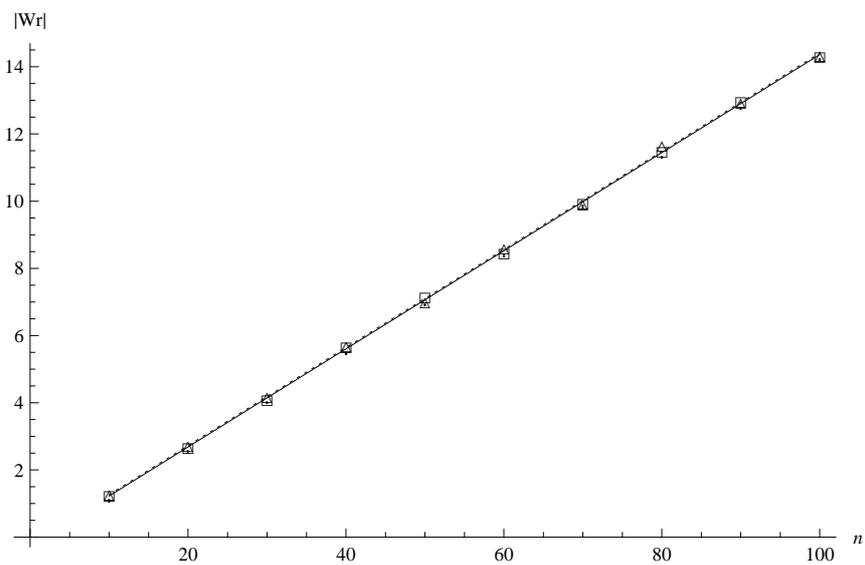}
     \caption{The mean absolute writhe of uniform random walks and polygons in confined space. Values obtained by computer simulations are shown by triangles and squares respectively. The black curve corresponds to the mean absolute value of the writhe of a uniform random polygon and the dashed curve corresponds to the mean absolute value of the writhe of a uniform random walk (open chain)}
     \label{urpmeanabswr}
   \end{center}
\end{figure}

We next estimate the mean absolute value of the writhe of an oriented uniform random walk or polygon.
Results are shown in Figure \ref{urpmeanabswr}.The black curve in the figure represents the mean absolute writhe of a uniform random polygon contained in $C^3$ with respect to its number of edges. The curve is fitted to a function of the form $a+bn$ where $a$ is estimated to be $-0.2372\pm0.04262$ and $b$ is estimated to be $0.1460\pm0.0007$, with a coefficient of determination $R^2=0.9998$.
The dashed curve in the figure represents the mean absolute writhe of a uniform random walk contained in $C^3$ with respect to its number of edges. The curve is fitted to a function of the form $a+bn$ where $a$ is estimated to be $-0.1958\pm0.0630$ and $b$ is estimated to be $0.1460\pm0.0010$, with a coefficient of determination $R^2=0.9996$.
Note that in both cases $b\approx\sqrt{q}$.

As the number of edges of the polygons increases, we observe a growth at a rate $O(n)$. This suggests the conjecture

\begin{con}
\begin{equation}
\sqrt{E[Wr^2]}\sim E[\sqrt{Wr^2}]=E[\left|Wr\right|]
\end{equation}
\end{con}

Its proof would strengthen our intuition that the mean writhe of a uniform random walk or polygon in confined space follows the normal distribution.

\bigskip

\subsubsection{Mean squared self-linking number and mean absolute self-linking number of an oriented uniform random walk or polygon in a confined space}

Our second numerical study concerns the self-linking of an oriented uniform random walk or polygon in confined space.

The self-linking number is defined as (Banchoff 1976),

\begin{equation}\label{sl2}
SL(l)=Wr(l)+\frac{1}{2\pi}\tau(l).
\end{equation}

The first term is the writhe of the curve, which we studied in the previous subsection and the second term is the total torsion of the curve, $\tau(l)$ divided by $2\pi$ which measures how much the curve deviates from being planar.

We calculated the mean squared self-linking number of an oriented uniform random walk of varying length and of an oriented uniform random polygon of varying length.

\begin{figure}
   \begin{center}
     \includegraphics{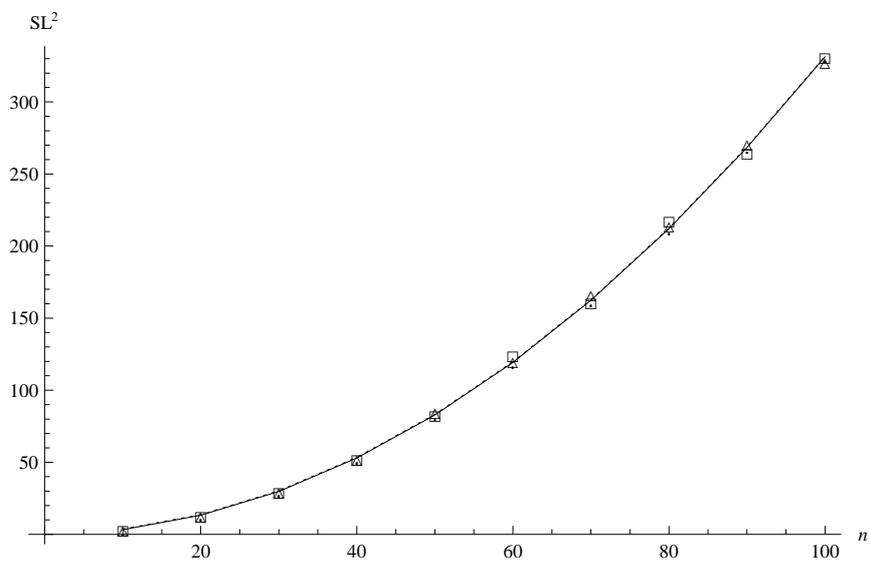}
     \caption{The mean squared self-linking number of uniform random walks and polygons. Values obtained by computer simulations are shown by triangles and squares respectively. The black curve is the graph of the mean squared self-linking number of uniform random polygons contained in $C^3$. The dashed curve is the graph of the mean squared self-linking number of uniform random walks contained in $C^3$.}
     \label{meansqsl}
   \end{center}
\end{figure}

Results are shown in Figure \ref{meansqsl}. The black curve in the figure illustrates the mean squared self-linking number of uniform random polygons in confined space and is fitted to a function of the form $qn^2+a$ where $q$ is estimated to be $0.0331\pm0.0003$ and $a$ is estimated to be $-0.0131\pm1.5226$, with a coefficient of determination $R^2=0.9993$.
The dashed curve in the figure illustrates the mean squared self-linking number of uniform random walks in confined space and is fitted to a function of the form $qn^2+a$ where $q$ is estimated to be $0.0331\pm0.0002$ and $a$ is estimated to be $0.3980\pm1.1778$, with a coefficient of determination $R^2=0.9996$. Thus the estimates given in the theorem are strongly supported by the data.

Note that the estimated value of $q$ coincides with the predicted value of $q$ given by Arsuaga \textit{et al} (2007), i.e. $q=0.0338\pm0.024$.

\begin{figure}
   \begin{center}
     \includegraphics{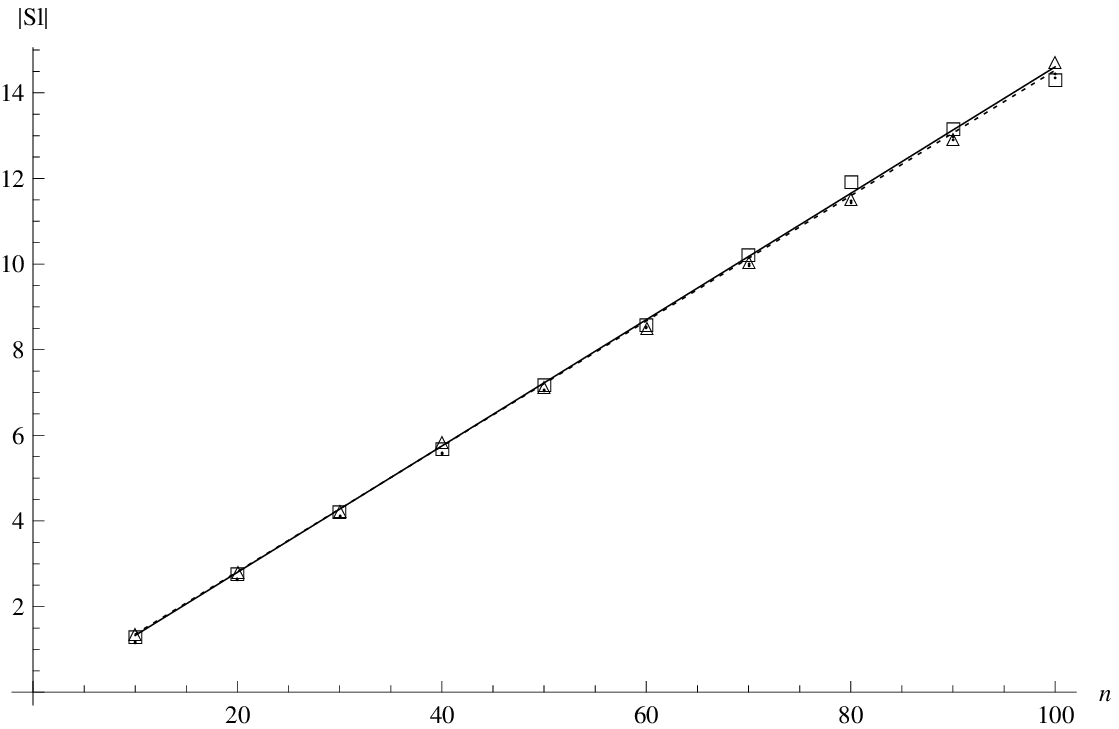}
     \caption{Mean absolute self-linking number of uniform random walks and polygons in confined space. Values obtained by computer simulations are shown by triangles and squares respectively. The black curve is the graph of the mean absolute self-linking number of uniform random polygons and the dashed curve is the graph of the mean absolute self-linking number of uniform random walks.}
     \label{abssl}
   \end{center}
\end{figure}

We next estimate the mean absolute value of the self-linking number of an oriented uniform random walk or polygon.
Results are shown in Figure \ref{abssl}. The black curve in the figure illustrates the mean absolute self-linking number of uniform random polygons contained in $C^3$ and is fitted to a curve of the form $a+bn$ where $a$ is estimated to be $-0.1516\pm0.1051$ and $b$ is estimated to be $0.1476\pm0.0017$, with a coefficient of determination $R^2=0.9990$.
The dashed curve in the figure illustrates the mean absolute self-linking number of uniform random walks contained in $C^3$ and is fitted to a curve of the form $a+bn$ where $a$ is estimated to be $-0.1089\pm0.0809$ and $b$ is estimated to be $0.1463\pm0.0013$, with a coefficient of determination $R^2=0.9994$. Note that in both cases $b\approx\sqrt{q}$.This suggests the conjecture:

\begin{con}
\begin{equation}
\sqrt{E[SL^2]}\sim E[\sqrt{SL^2}]=E[\left|SL\right|]
\end{equation}
\end{con}

Its proof would suggest that the mean self-linking number of a uniform random walk or polygon follows the normal distribution.

\subsubsection{Mean squared and mean absolute linking number of oriented uniform random walks in a confined space}

Our third numerical study concerns the linking between two oriented uniform random walks or polygons in confined space.
By Theorem \ref{uniflink}, the mean squared linking number between two oriented uniform random polygons grows at a rate $O(n^2)$ with regard to the number of edges of the polygons.

\begin{figure}
   \begin{center}
     \includegraphics{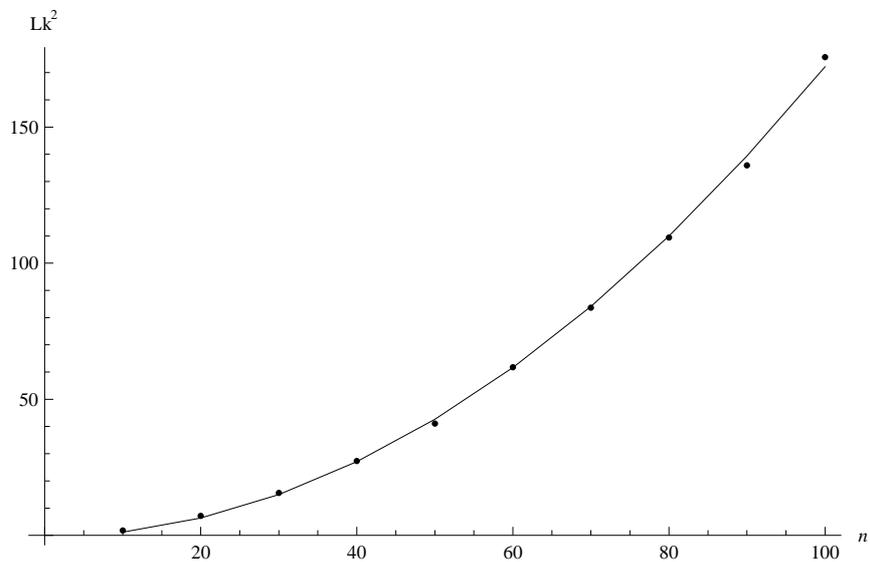}
     \caption{The mean squared linking number of two uniform random walks contained in $C^3$.}
     \label{f}
   \end{center}
\end{figure}

We calculated the mean squared linking number between two oriented uniform random walks of varying lengths and the results are shown in Figure \ref{f}. The curve in the figure is fitted to a function of the form $a+\frac{q}{2}n^2$ where $a$ is estimated to be $-0.5499\pm0.9297$ and $q$ is estimated to be $0.0346\pm0.0004$, with a coefficient of determination $R^2=0.9991$ . Clearly, as the number of edges of the walks increases, we observe a growth at a rate $O(n^2)$, as expected from Theorem \ref{uniflink2}.

\begin{figure}
   \begin{center}
     \includegraphics{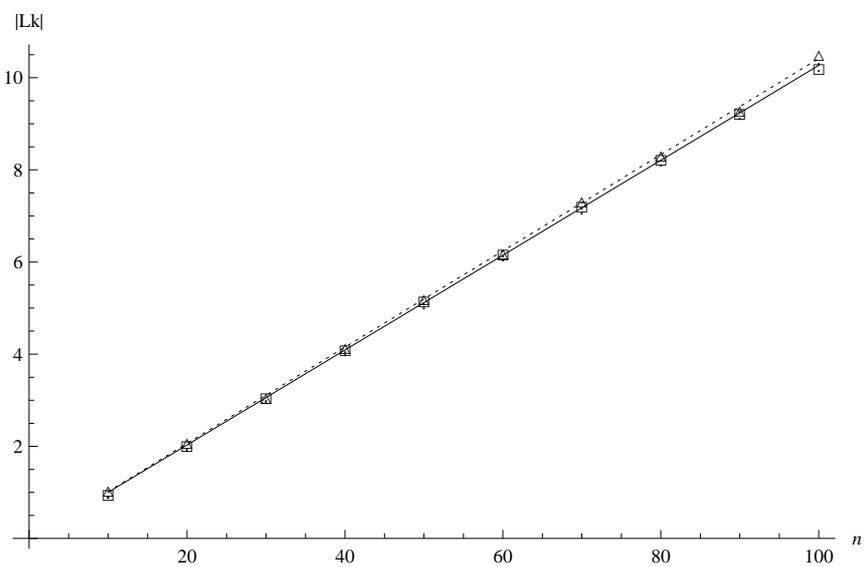}
     \caption{The mean absolute linking number of two uniform random walks and polygons contained in $C^3$. Values obtained by computer simulations are shown by triangles and squares respectively. The black curve is the graph of the mean absolute linking number of two uniform random polygons and the dashed curve is the graph of the mean absolute linking number of two uniform random walks.}
     \label{abslkunif}
   \end{center}
\end{figure}
\bigskip

Next we calculated the mean absolute linking number between two oriented uniform random walks or polygons of varying length and results are shown in Figure \ref{abslkunif}. The black curve in the figure illustrates the mean absolute linking number of two uniform random polygons contained in $C^3$ and is fitted to a function of the form $a+bn$ where $a$ is estimated to be $-0.0300\pm0.0256$ and $b$ is estimated to be $0.1030\pm0.0004$, with a coefficient of determination $R^2=0.9999$ . The dashed curve in the figure illustrates the mean absolute linking number of two uniform random walks contained in $C^3$ and is fitted to a function of the form $a+bn$ where $a$ is estimated to be $-0.0251\pm0.0350$ and $b$ is estimated to be $0.1044\pm0.0006$, with a coefficient of determination $R^2=0.9998$ . Clearly, as the number of edges of the polygons increases, we observe a growth at a rate $O(n)$. This suggests the conjecture:

\begin{con}
\begin{equation}\label{mean}
\sqrt{E[Lk^2]}\sim E[\sqrt{Lk^2}]=E[\left|Lk\right|].
\end{equation}
\end{con}

\begin{rem}\rm Note that (\ref{mean}) agrees with our intuition and with the numerical results obtained by Arsuaga \textit{et al} (2007) that the linking number of two uniform random walks or polygons in confined space follows the normal distribution. Indeed, if $Lk\sim N(0,\sigma^2)$, where $\sigma^2=O(n^2)$, then $|Lk|$ follows the half-normal distribution and $E[|Lk|]=\sigma\sqrt{2/\pi}=O(n)$. Then $\frac{|Lk|}{\sigma}$ follows the $\chi$-distribution and $Lk^2/\sigma^2$ follows the $\chi^2$-distribution with one degree of freedom and mean $1$. Thus, $E[Lk^2/\sigma^2]=1$, $E[Lk^2]=\sigma^2=O(n^2)$, hence we have that $E[Lk^2]\approx E[|Lk|]^2$.
\end{rem}

\begin{figure}
   \begin{center}
     \includegraphics{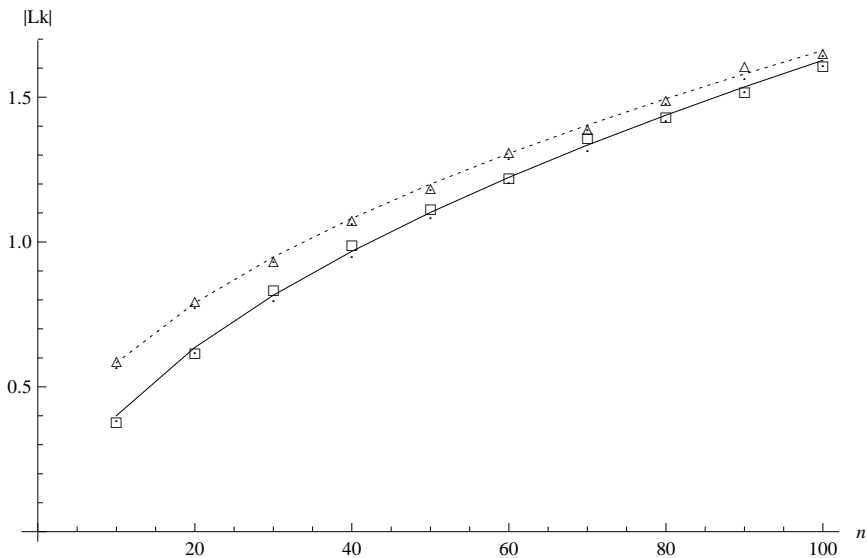}
     \caption{Mean absolute linking number of a uniform random walk or polygon and a fixed square in confined space. Values obtained by computer simulations are shown by triangles and squares respectively. The black curve is the graph of the mean absolute linking number of a uniform random polygon and a fixed square and the dashed curve is the graph of the mean absolute linking number of a uniform random walk and a fixed square.}
     \label{circleuniflink}
   \end{center}
\end{figure}

Our next numerical study concerns the linking between an oriented uniform random walk or polygon and a fixed simple closed curve in confined space.
By Theorem \ref{abs2} the mean absolute value of the linking number between an oriented uniform walk or polygon of $n$ edges and a fixed oriented simple closed curve in confined space has a scaling of the form $O(\sqrt{n})$. First we consider the oriented square $S_1$ defined by the sequence of vertices $(0.1,0.1,0.5),(0.9,0.1,0.5),(0.9,0.9,0.5)$,
$(0.1,0.9,0.5),(0.1,0.1,0.5)$ and a uniform random walk or polygon. The results of our simulations can be seen in Figure \ref{circleuniflink}. The black curve shows the growth rate of the mean absolute linking number of a uniform random polygon and $S_1$. The data are fitted to a function of the form $a+b\sqrt{n}$, where $a$ is estimated to be $-0.1665\pm0.0212$ and $b$ is estimated to be $0.1794\pm 0.0029$, with a coefficient of determination $R^2=0.9980$. The dashed curve shows the growth rate of the mean absolute linking number of a uniform random walk and $S_1$.
The data are fitted to a function of the form $a+b\sqrt{n}$, where $a$ is estimated to be $0.0848\pm0.0153$ and $b$ is estimated to be $0.1576\pm0.0021$, with a coefficient of determination $R^2=0.9986$. Thus, we can see that the data is consistent with Theorem \ref{abs2}.

\begin{figure}
   \begin{center}
     \includegraphics{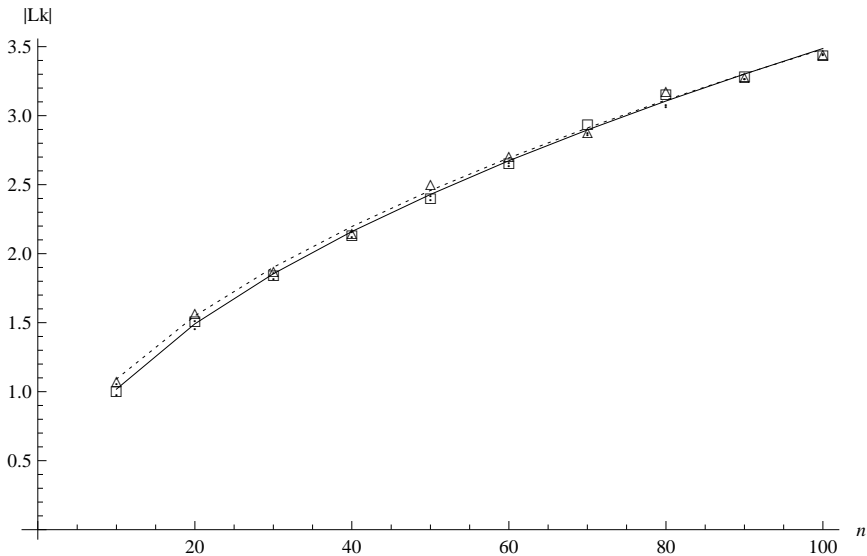}
     \caption{The mean absolute linking number of a fixed trefoil and a uniform random walk or polygon in confined space. Values obtained by computer simulations are shown by triangles and squares respectively. The black curve is the graph of the mean absolute linking number of a fixed trefoil and a uniform random polygon. The dashed curve is the graph of the mean absolute linking number of a fixed trefoil and a uniform random walk.}
     \label{circleuniflink2}
   \end{center}
\end{figure}

To illustrate that the rate of growth of the mean absolute linking number does not depend on the knot type of the fixed simple closed curve, we consider the oriented trefoil $S_2$ defined by the sequence of vertices $(0.9,0.5,0.5),(0.1,0.5,0.4)$, $(0.5,0.3,0.9),(0.6,0.3,0.1),(0.2,0.9,0.6), (0.5,0.2,0.5)$, $(0.9,0.5,0.5)$ and a uniform random walk or polygon. The results of our simulations can be seen in Figure \ref{circleuniflink2}. The black curve shows the growth rate of the mean absolute linking number of a uniform random polygon and $S_2$. The data are fitted to a function of the form $a+b\sqrt{n}$, where $a$ is estimated to be $-0.1232\pm0.0381$ and $b$ is estimated to be $0.3610\pm0.00514$, with a coefficient of determination $R^2=0.9984$. The dashed curve shows the growth rate of the mean absolute linking number of a uniform random walk and $S_2$.
The data are fitted to a function of the form $a+b\sqrt{n}$, where $a$ is estimated to be $-0.0142\pm0.0432$ and $b$ is estimated to be $0.3496\pm0.0058$, with a coefficient of determination $R^2=0.9978$. Hence, the data confirm our analytical result in Theorem \ref{abs2} for any fixed simple closed curve.

\bigskip

\subsection{Numerical Results on Equilateral Random Walks}
Equilateral random walks are widely used to study the behavior of polymers under $\theta$-conditions. It is of great interest to study the scaling of the linking number, the self-linking number and the writhe of equilateral random walks and polygons. These situations are also much more complex in comparison to the uniform models studied in previous sections, since the probability of crossing of two edges depends on their distance and the probability of positive or negative crossing is independent upon the previous edges. In this section we present numerical results concerning equilateral random walks and polygons. It would be of great interest to have rigorous proofs of the scaling observed.

\subsubsection{Mean absolute self-linking number of an equilateral random walk}

\bigskip

\begin{figure}
   \begin{center}
     \includegraphics{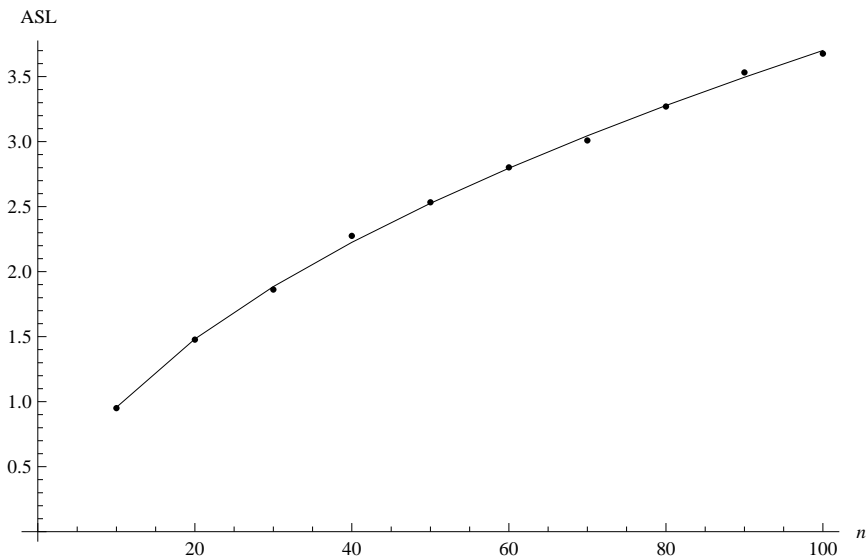}
     \caption{Mean absolute value of the self-linking number of an equilateral random walk}
     \label{rwmeanabswr}
   \end{center}
\end{figure}

\bigskip

In this section we discuss our numerical results on the scaling of the self-linking number of a random walk.

The self-linking number of a random walk $X$ is equal to the writhe plus the total torsion of the random walk.
\begin{equation}
\begin{split}
&E[|SL(X)|]=E[|Wr(X)+\tau(X)|]\leq E[|Wr(X)|]+E[|\tau(X)|]\\
&=E[|Wr(X)|]+E[|\sum_{i}\phi_i|]\leq E[|Wr(X)|]+E[\sum_{i}|\phi_i|]
\end{split}
\end{equation}

where $E[\sum_{i}|\phi_i|]$ has been proved to be approximately equal to $n\frac{\pi}{2}-\frac{3\pi}{8}$ (Plunkett \textit{et al} 2007, Grosberg 2008).

Previous numerical results (Orlandini \textit{et al} 1994) suggest that the average of the absolute value of the writhe of an ideal walk increases as $\sqrt{n}$, where $n$ is the length of the walk.

Figure \ref{rwmeanabswr} shows the $\langle\left|ASL\right|\rangle$ values obtained in numerical simulations of ideal random walks.

We have fitted the computing data points with the function $a+b\sqrt{n}$, leaving the two parameters $a$ and $b$ free. Then $a$ was estimated to be $-0.3131\pm0.0313$ and $b$ was estimated to be $0.4014\pm0.0042$, with a coefficient of determination $R^2=0.9991$.

\bigskip

\subsubsection{Mean absolute linking number of two equilateral random walks whose starting points coincide}

\bigskip

\begin{figure}
   \begin{center}
     \includegraphics{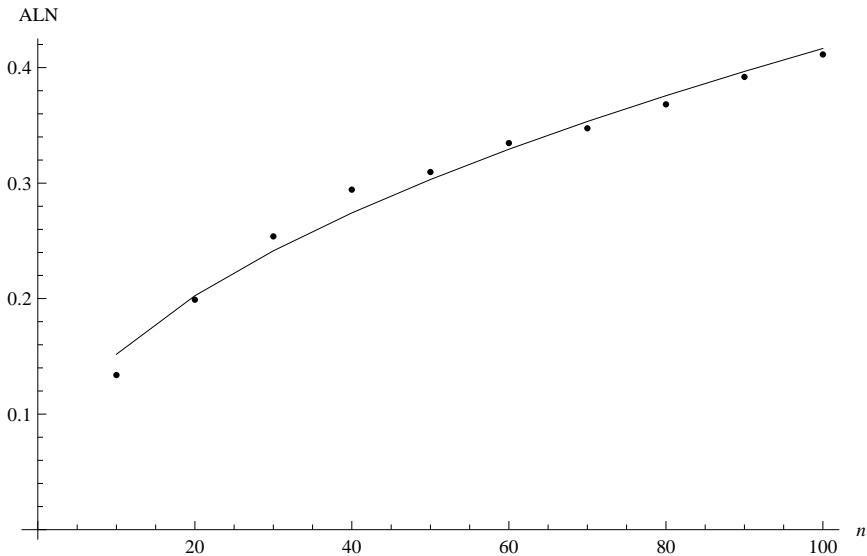}
     \caption{The mean absolute value of the linking number of two equilateral random walks whose starting points coincide}
     \label{rwmeanabslk}
   \end{center}
\end{figure}

In this section we discuss our numerical results on the scaling of the linking number between two equilateral random walks whose starting points coincide.
Figure \ref{rwmeanabslk} shows the $\langle ALN\rangle$ values obtained in numerical simulations of ideal random walks in a non-constrained linear form.

By the numerical results presented in the last section, we expect that the mean absolute value of the writhe of an oriented equilateral random walk $X$ of $n$ steps will have a power-law dependence on the length of the walk:

\begin{equation}
\langle AWR\rangle\approx O(\sqrt{n}).
\end{equation}

\bigskip

Let $X=(X_0,X_1,\dots,X_n)$ and $Y=(Y_0,Y_1,\dots,Y_n)$ denote two oriented equilateral random walks of length $n$ whose starting points coincide, i.e. $X_0=Y_0=0$. One can use the scaling of $\langle AWR\rangle$ to give information concerning the scaling of $\langle ALN\rangle$ as follows :

We define $X-Y$ to be the oriented equilateral random walk of $2n$ steps $(Y_n,\dots,Y_1,Y_0=X_0,X_1,\dots,X_n)$. Its writhe then is

\begin{eqnarray}
&Wr(X-Y)=Wr(X)+Wr(-Y)+2L(X,-Y) \nonumber\\
&\Leftrightarrow Wr(X-Y)=Wr(X)+Wr(-Y)-2L(X,Y) \nonumber\\
&\Rightarrow 2L(X,Y)=-Wr(X-Y)+Wr(X)+Wr(Y) \nonumber\\
&\Rightarrow\left|L(X,Y)\right|\leq\frac{1}{2}\bigl(\left|Wr(X-Y)\right|+\left|Wr(X)\right|+\left|Wr(Y)\right|\bigr)\nonumber\\
&\Rightarrow\left|L(X,Y)\right|\le  \frac{1}{2}\bigl(O(\sqrt{2n})+ O(\sqrt{n})+O(\sqrt{n})\bigr)\nonumber\\
&\Rightarrow\left|L(X,Y)\right|\le O(\sqrt{n}).
\end{eqnarray}

We decided therefore to check whether the average of the absolute value of the linking number between two ideal walks increases as $\sqrt{n}$. We have fitted the computing data points with the function $a+b\sqrt{n}$, leaving the two parameters $a$ and $b$ free. Then $a$ was estimated to be $0.0294\pm0.0130$ and $b$ was estimated to be $0.0387\pm0.0018$ with a coefficient of determination $0.9839$.

\bigskip

\section{Conclusions}

\bigskip

The measurement of the entanglement of open chains is of great interest for many applications, such as the study of the properties of polymer melts. In this paper, we focused our study in the case of uniform random walks (open chains) and polygons in confined volumes. In Theorems \ref{unifsl}, \ref{uniflink}, \ref{unisl} we gave rigorous proofs that the scaling of the mean squared linking number, the mean squared writhe and the mean squared self-linking number of oriented uniform random walks and polygons in confined space, with respect to their length, is of the form $O(n^2)$.

Further, we are interested in the mean absolute value of the linking number of two uniform random walks or polygons in confined space. In this direction, we prove in Theorem \ref{abs2} that the mean absolute value of the linking number of an oriented uniform random walk or polygon and a fixed oriented simple closed curve in confined space is of the form $O(\sqrt{n})$. Our numerical results confirm the analytical prediction, and furthermore suggest that for two oriented uniform random walks or polygons in confined space $\sqrt{E[lk^2]}\sim O(n)\sim E[\sqrt{lk^2}]$. A possible direction for future work would be to prove these results analytically.

Ideal random walks are used to model the behaviour of polymers under $\theta$-conditions. We have analyzed numerically the scaling of the mean absolute value of the linking number between two equilateral random walks of $n$ steps and the mean absolute value of the self-linking number of an equilateral random walk of $n$ steps. Both appear to scale as $O(\sqrt{n})$. An important direction for future work is to complete the analysis and provide proofs for the scaling of the self-linking number and that of the linking number of equilateral random walks and polygons.

\end{document}